\documentclass[11pt]{article}

\usepackage{epsfig}
\usepackage{graphicx}
\usepackage{amsbsy}
\usepackage{amsmath}
\usepackage{amsfonts}
\usepackage{amssymb}
\usepackage{textcomp}
\usepackage{hyperref}
\usepackage{aliascnt}

\newcommand{\mcm}[3]{\newcommand{#1}[#2]{{\ensuremath{#3}}}} 

\mcm{\tuple}{1}{\langle #1 \rangle}
\mcm{\name}{1}{\ulcorner #1 \urcorner}
\mcm{\Nbb}{0}{\mathbb{N}}
\mcm{\Zbb}{0}{\mathbb{Z}}
\mcm{\Rbb}{0}{\mathbb{R}}
\mcm{\Cbb}{0}{\mathbb{C}}
\mcm{\Bcal}{0}{\cal B}
\mcm{\Ccal}{0}{\cal C}
\mcm{\Dcal}{0}{\cal D}
\mcm{\Ecal}{0}{\cal E}
\mcm{\Fcal}{0}{\cal F}
\mcm{\Gcal}{0}{\cal G}
\mcm{\Hcal}{0}{\cal H}
\mcm{\Ical}{0}{\cal I}
\mcm{\Lcal}{0}{\cal L}
\mcm{\Mcal}{0}{\cal M}
\mcm{\Ncal}{0}{\cal N}
\mcm{\Pcal}{0}{{\cal P}}
\mcm{\Scal}{0}{{\cal S}}
\mcm{\Tcal}{0}{{\cal T}}
\mcm{\Ucal}{0}{{\cal U}}
\mcm{\Vcal}{0}{{\cal V}}
\mcm{\Ycal}{0}{{\cal Y}}
\mcm{\Mfrak}{0}{\mathfrak M}
\newcommand{\easy}{\nobreak\hfill$\square$\medskip}

\DeclareMathOperator{\Cl}{Cl}
\mcm{\restric}{0}{\upharpoonright}
\mcm{\upset}{0}{\uparrow}
\mcm{\onto}{0}{\twoheadrightarrow}
\mcm{\smallNbb}{0}{{\small \mathbb{N}}}
\DeclareMathOperator{\preop}{op}
\mcm{\op}{0}{^{\preop}}
\newcommand{\supth}{\textsuperscript{th}}

{\begin{array}{c}
\setlength{\unitlength}{1em}}%
{\end{array}}

\usepackage{amsthm}

\newcommand{\theoremize}[2]{\newaliascnt{#1}{thm} \newtheorem{#1}[#1]{#2} \aliascntresetthe{#1}}

\theoremstyle{plain}
\newtheorem{thm}{Theorem}[section]
\theoremize{lem}{Lemma}
\theoremize{sublem}{Sublemma}
\theoremize{claim}{Claim}
\theoremize{rem}{Remark}
\theoremize{prop}{Proposition}
\theoremize{cor}{Corollary}
\theoremize{que}{Question}
\theoremize{oque}{Open Question}
\theoremize{con}{Conjecture}

\theoremstyle{definition}
\theoremize{dfn}{Definition}
\theoremize{eg}{Example}
\theoremize{exercise}{Exercise}
\theoremstyle{plain}

\usepackage{verbatim}
\usepackage[all]{xy}

\usepackage{subfig}

\title{\scshape %A connection between 
Matroid intersection, \\ base packing and base covering\\ for infinite matroids}
\author{Nathan Bowler \and Johannes Carmesin }

\begin{document}

%TODO at Nathan, Do you know a better title

\maketitle
 \begin{abstract}
As part of the recent developments in infinite matroid theory, there have been a number of conjectures about how standard theorems of finite matroid theory might extend to the infinite setting. These include base packing, base covering, and matroid intersection and union. We show that several of these conjectures are equivalent, so that each gives a perspective on the same central problem of infinite matroid theory.
For finite matroids, these equivalences give new and simpler proofs for the finite theorems corresponding to these conjectures.

This new point of view also allows us to extend, and simplify the proofs of, some cases where these conjectures were known to be true.
 \end{abstract}

\newcommand{\cw}{at most countably weird}

\section{Introduction}

The well-known finite matroid intersection theorem of Edmonds states that for any two finite matroids $M$ and $N$ the size
of a biggest common independent set is equal to the minimum of the rank sum $r_{M}(E_M)+r_{N}(E_N)$, where the minimum is taken over all partitions 
 $E=E_M \dot\cup E_N$. The same statement for infinite matroids is true, but for a silly reason \cite{christian_phd}, which suggests that more care is needed in extending this statement to the infinite case. 

Nash-Williams \cite{Aharoni:Ziv:92} proposed the following for finitary matroids.

\begin{con}\label{int}
Any two matroids $M$ and $N$ on a 
common ground set $E$ have a common independent set $I$ admitting a partition $I
= J_M \cup J_N$ such that $cl_{M}(J_M) \cup cl_{N}(J_N) = E$.
\end{con}

 For finite matroids this is easily seen to be equivalent to the intersection theorem, so we will refer to \autoref{int} as the Matroid Intersection Conjecture.
 In \cite{union2}, it was shown that this conjecture implies the celebrated Aharoni-Berger-Theorem \cite{AharoniBerger}, also known as the Erd\H{o}s-Menger-Conjecture.
Call a matroid \emph{finitary} if all its circuits are finite and \emph{co-finitary} if its dual is finitary.
The conjecture is true in the cases where $M$ is finitary and $N$ is co-finitary \cite{union2}.%
\footnote{In fact in \cite{union2} the conjecture was proved for a slightly larger class.}
Aharoni and Ziv \cite{Aharoni:Ziv:92} proved the conjecture for one matroid finitary and the other a countable direct
sum of finite rank matroids.

In this paper we will demonstrate that the Matroid Intersection Conjecture is a natural formulation by showing that it is equivalent 
to several other new conjectures in unexpectedly different parts of infinite matroid theory.

Suppose we have a family of matroids $(M_k | k \in K)$ on the same ground set $E$.
A  {\em packing} for this family consists of a spanning set $S_k$ for each $M_k$ such that the $S_k$ are all disjoint.
Note that not all families of matroids have a packing.
More precisely, the well-known finite base packing theorem states that if $E$ is finite then the family has a packing if and only if for every subset $Y\subseteq E$ the following holds.
\[
\sum_{k\in K}  r_{M_k.Y}(Y) \leq |Y|
\]
The Aharoni-Thomassen-graphs \cite{aharoniThom, DiestelBook10} show that this theorem does not extend verbatim to finitary matroids. 
However, the base packing theorem extends to finite families of co-finitary matroids \cite{union1}. This implies the topological tree packing theorems of Diestel and Tutte.
Independently from our main result, we close the gap in between by showing that the base packing theorem extends to arbitrary families of co-finitary matroids (for example, topological cycle matroids).

Similar to packings are coverings: a {\em covering} for the family $(M_k | k \in K)$ consists of an independent set $I_k$ for each $M_k$ such that the $I_k$ cover $E$. 
And analoguesly to the base packing theorem, there is a base covering theorem characterising the finite families of finite matroids admitting a covering.

We are now in a position to state our main conjecture, which we will show is equivalent to the intersection conjecture.
Roughly, the finite base packing theorem says that a family has a packing if it is very dense. Similarly, 
the finite base covering theorem says roughly that a family has a covering if it is very sparse.
Although not every family of matroids has a packing and not every family has a covering, we could ask if it is always possible to
divide the ground set into a ``dense'' part, which has a packing, and a ``sparse'' part, which has a covering?
 More precisely, we conjecture the following:

\begin{con}\label{pc}
For any family of matroids $(M_k | k \in K)$ on the same ground set $E$, the ground set admits a partition $E=P\dot\cup C$
such that $(M_k\restric_P | k \in K)$  has a packing and $(M_k.C | k \in K)$ has a covering.
\end{con}

Here $M_k\restric_P$ is the restriction of $M_k$ to $P$ and $M_k.C$ is the contraction of $M_k$ onto $C$.
Note that if $(M_k\restric_P | k \in K)$ has a packing, then $(M_k.P | k \in K)$ has a packing, so we get a stronger statement by taking the restriction here. Similarly, we get a stronger statement by contracting to get the family which should have a covering than we would get by restricting.

For finite matroids, we show that this new conjecture is true and implies the base packing and base covering theorems.
So the finite version of \autoref{pc} unifies the base packing and the base covering theorem into one theorem.

For infinite matroids, we show that \autoref{pc} and the intersection conjecture are equivalent, and that both are equivalent to \autoref{pc} for pairs of matroids.
As \autoref{pc} for pairs of matroids is self-dual, this shows the less obvious fact that the intersection conjecture is self-dual:

\begin{cor}\label{int_dual}
If $M$ and $N$ are matroids on the same ground set then \autoref{int} is true for $M$ and $N$ iff it is true for $M^*$ and $N^*$.
\end{cor}

\autoref{pc} also suggests a base packing conjecture and a base covering conjecture which we show are equivalent 
to the intersection conjecture but not to the above mentioned rank formula formulation of base packing for infinite matroids.

The various results about when intersection is true transfer via these equivalences to give results 
showing that these new conjectures also hold in the corresponding special cases.
For example, while the rank-formulation of the covering theorem is not true for all families of cofinitary matroids, the new covering conjecture is true  in that case.
This yields a base covering theorem for the algebraic cycle matroid of any locally finite graph and the topological cycle matroid of any graph. Similarly, we immediately obtain in this way that the new packing and covering conjectures are true for finite families of finitary matroids. Thus we get packing and covering theorems for the finite cycle matroid of any graph.

For finite matroids, the proofs of the equivalences of these conjectures simplify the proofs of the corresponding finite theorems.

We show that \autoref{pc} might be seen as the infinite analogue of the rank formula of the matroid union theorem.
It should be noted that there are two matroids whose union is not a matroid \cite{union1}, so there is no infinite analogue of the finite matroid union theorem as a whole.

This new point of view also allows us to give a simplified account of the special cases of the intersection conjecture and even to extend the results a little bit. Our result includes the following:

\begin{thm}
For any family of matroids $(M_k | k \in K)$ on the same ground set $E$ which between them have only countably many circuits, the ground set admits a partition $E=P\dot\cup C$
such that $(M_k\restric_P | k \in K)$  has a packing and $(M_k.C | k \in K)$ has a covering.
\end{thm}

This paper is organised as follows:
In Section 2, we recall some basic matroid theory and introduce a key idea, that of exchange chains.
After this, in Section 3, we restate our main conjecture and look at its relation to the infinite matroid intersection conjecture.
In Section 4, we prove a special case of our main conjecture.
In the next two sections, we consider base coverings and base packings of infinite matroids.
In the final section, Section 7, we give an overview over the various equivalences we have proved.

\section{Preliminaries}

\subsection{Basic matroid theory}

Throughout, notation and terminology for graphs are that of~\cite{DiestelBook10}, for matroids  that of~\cite{Oxley,matroid_axioms}, and for topology that of~\cite{Armstrong}.
$M$ always denotes a matroid and $E(M)$, $\Ical(M)$, $\Bcal(M)$, $\Ccal(M)$ and $\Scal(M)$ denote its ground set and its sets of independent sets, bases, circuits and spanning sets, respectively.

Recall that the set $\Ical(M)$ is required to satisfy the following\emph{independence axioms}~\cite{matroid_axioms}:
\begin{itemize}
	\item[(I1)] $\emptyset\in \Ical(M)$.
	\item[(I2)] $\Ical(M)$ is closed under taking subsets.
	\item[(I3)] Whenever $I,I'\in \Ical(M)$ with $I'$ maximal and $I$ not maximal, there exists an $x\in I'\setminus I$ such that $I+x\in \Ical(M)$.
	\item[(IM)] Whenever $I\subseteq X\subseteq E$ and $I\in\Ical(M)$, the set $\{I'\in\Ical(M)\mid I\subseteq I'\subseteq X\}$ has a maximal element.
\end{itemize}

The axiom (IM) for the dual $M^*$ of $M$ is equivalent to the following: 
\begin{itemize}
		\item[(IM$^*$)] Whenever $Y\subseteq S\subseteq E$ and $S\in \Scal(M)$, the set $\{S'\in\Scal(M)\mid Y\subseteq S'\subseteq S\}$ has a minimal element.
\end{itemize}
As the dual of any matroid is also a matroid, every matroid satisfies this.
We need the following facts about circuits \cite{matroid_axioms}:

\begin{enumerate}
	\item [(C3)] Whenever $X\subseteq C\in \Ccal(M)$ and $\{C_x\mid x\in X\} \subseteq \Ccal(M)$ satisfies $x\in C_y\Leftrightarrow x=y$ for all $x,y\in X$, 
then for every $z \in C\setminus \left( \bigcup_{x \in X} C_x\right)$ there exists a  $C'\in \Ccal(M)$ such that $z\in C'\subseteq \left(C\cup  \bigcup_{x \in X} C_x\right) \setminus X$.
        \item [(C4)] Every dependent set contains a circuit. \label{exists_cir}
\end{enumerate}

A matroid is called \emph{finitary} if every circuit is finite.
An \emph{$M$-bond} is a circuit of $M^*$.

\begin{lem}\label{meetbond}
A set $S$ is $M$-spanning iff it meets every $M$-bond.
\end{lem}

\begin{proof}
We prove the dual version where $I:=E(M)\setminus S$.
\begin{equation}\label{ind_version}
\begin{minipage}[c]{0.8\textwidth}
 A set $I$ is $M^*$-independent iff it does not contain an $M^*$-circuit.
\end{minipage}
\end{equation}
Clearly, if $I$ contains a circuit, then it is not independent.
Conversely, if $I$ is not independent, then by (C4) it also contains a circuit.
\end{proof}

Let $2^X$ denote the power set of $X$.
If $M=(E,\Ical)$ is a matroid, then for every $X\subseteq E$ the \emph{restriction matroid $M\restric_X:=(X,\Ical\cap 2^X)$},
the \emph{deletion matroid $M-X:=M\restric_{E-X}$},
the \emph{contraction matroid $M.X:=(M^*\restric_X)^* $} and
 the \emph{contracted matroid $M/X:=M.(E-X)$},
are also matroids.

\begin{lem}\label{sthg}
 Let $M$ be a matroid and $X\subseteq E(M)$. If $S_1\subseteq X$ spans $M\restric_X$ and $S_2\subseteq E\setminus X$ spans $M/X$, then $S_1\cup S_2$ spans $M$.
\end{lem}

\begin{proof}
We will apply \autoref{meetbond}: so let $B$ be any $M$-bond.
If $B\cap X$ is nonempty, then it is easy to see that $B\cap X$ contains an $M\restric_X$-bond, so $S_1$ meets $B$.
Otherwise $B\subseteq E-X$, and it suffices to show that $S_2$ meets $B$, that is $B$ is an $M/X$-bond.
This follows from the fact that $B$ is an $M^*$-circuit, so also an $M^*\restric_{(E-X)}$-circuit.
\end{proof}

\begin{lem}[\cite{bruhn:wollan_con}, Lemma]\label{cir_cocir}
Let $M$ be a matroid with a circuit $C$ and a co-circuit $D$, then $|C\cap D|\neq 1$.
\end{lem}

A particular class of matroids we shall employ is the {\em uniform} matroids $U_{n,E}$ on a groundset $E$, in which the bases are the subsets of $E$ of size $n$. In fact, the matroids we will use are those of the form $U^*_{1,E}$, in which the bases are all those sets obtained by removing a single element from $E$. Such a matroid is said to consist of a single circuit, because $\Ccal(U^*_{1,E}) = \{E\}$. A subset is independent iff it isn't the whole of $E$. Note that for a subset $X$ of $E$, $U^*_{1,E} \restric_X$ is free (every subset is independent) unless $X$ is the whole of $E$, and $U^*_{1,E}.X = U^*_{1,X}$ unless $X$ is empty.

\subsection{Exchange chains}

\begin{figure}[htbp]
	\centering
	\subfloat[Before the exchange\label{fig:chain1}]{\includegraphics{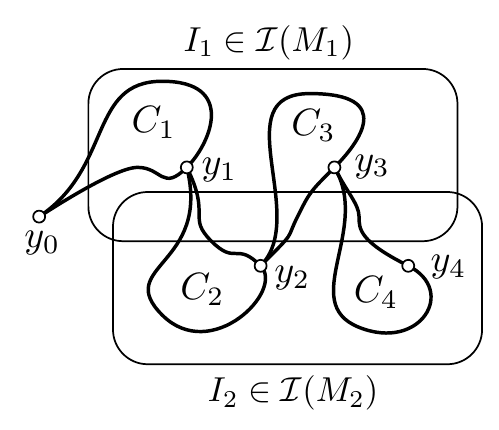}}
	\hspace{1cm}
	\subfloat[After the exchange \label{fig:chain2}]{\includegraphics{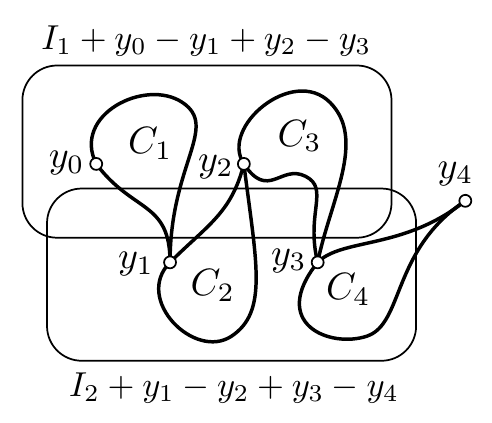}}
	\caption{An $(I_1,I_2)$-exchange chain of length $4$.}
	\label{fig:chain}
\end{figure}

Below, we will need a modification of the concept of exchange chains introduced in \cite{union1}.
The only modification is that we need not only exchange chains for families with two members but more generally exchange chains for arbitrary families,
which we define as follows:
Let $(M_k|k\in K)$ be a family of matroids and let $B_k\in \Ical(M_k)$.
A \emph{$(B_k|k\in K)$-exchange chain (from $y_0$ to $y_n$)} is a tuple $(y_0, k_0;y_1,k_1;\ldots;y_n)$ where
$B_{k_l}+y_l$ includes an $M_{k_l}$-circuit containing $y_l$ and $y_{l+1}$.
A $(B_k|k\in K)$-exchange chain from $y_0$ to $y_n$ is called \emph{shortest} 
if there is no $(B_k|k\in K)$-exchange chain $(y_0', k_0';y_1',k_1';\ldots;y_m')$ with $y_0'=y_0$, $y_m'=y_n$ and $m<n$.
A typical exchange chain is shown in \autoref{fig:chain}.

\begin{comment}
 \begin{lem}\label{runchains}
Let $(M_k|k\in K)$ be a family of matroids and let $B_k\in \Ical(M_k)$.
If $(y_0, k_0;y_1,k_1;\ldots;y_n)$ is a shortest $(B_k|k\in K)$-exchange chain from $y_0$ to $y_n$,
then $B_k'\in \Ical(M_k)$ for every $k$, where
\[
 B_k':=B_k\cup \{y_l|k_{l}=k\} \setminus \{y_{l+1}|k_{l}=k\}
\]
If $B_k$ spans $B_k'$, then $cl_{M_k}B_k=cl_{M_k}B_k'$.
\end{lem}

We only sketch the proof roughly.
The proof that the $B_k'$ are independent is done by induction on $n$ and is that of Lemma 4.2 in \cite{union1}.
To see the second assertion, we apply the following basic Lemma \cite{union2} to $M_k\restric_{cl_{M_k}}(B_k)$.
\begin{lem}\label{thm:sym-dif}
	Let $M$ be a matroid and $I, B\in\Ical(M)$ with $B$ maximal and $B\setminus I$ finite.
	Then $|I\setminus B| \leq |B\setminus I|$.
\end{lem} 
Note that $|B_k\setminus B_k'|$ is finite and is equal to $|B_k'\setminus B_k|$.
\end{comment}

 \begin{lem}\label{runchains}
Let $(M_k|k\in K)$ be a family of matroids and let $B_k\in \Ical(M_k)$.
If $(y_0, k_0;y_1,k_1;\ldots;y_n)$ is a shortest $(B_k|k\in K)$-exchange chain from $y_0$ to $y_n$,
then $B_k'\in \Ical(M_k)$ for every $k$, where
\[
 B_k':=B_k\cup \{y_l|k_{l}=k\} \setminus \{y_{l+1}|k_{l}=k\}
\]
Moreover, $cl_{M_k}B_k=cl_{M_k}B_k'$.
\end{lem}

\begin{proof}[Proof (Sketch)]
The proof that the $B_k'$ are independent is done by induction on $n$ and is that of Lemma 4.2 in \cite{union1}.
To see the second assertion, first note that $\{y_l|k_{l}=k\}\subseteq cl_{M_k}B_k$ and thus $B_k'\subseteq cl_{M_k}B_k$.
Thus it suffices to show that $B_k\subseteq cl_{M_k}B_k'$. To see this, note that $|B_k\setminus B_k'|$ is finite and is equal to $|B_k'\setminus B_k|$
and conclude that $B_k'$ is a base of $M_k\restric_{cl_{M_k}(B_k)}$ by the following basic Lemma \cite{union2} applied with $M=M_k\restric_{cl_{M_k}(B_k)}$, $B=B_k$, 
$I$ a base of $M_k\restric_{cl_{M_k}(B_k)}$ containing $B_k'$.
\end{proof}

\begin{lem}\label{thm:sym-dif}
	Let $M$ be a matroid and $I, B\in\Ical(M)$ with $B$ maximal and $B\setminus I$ finite.
	Then $|I\setminus B| \leq |B\setminus I|$.
\end{lem}

\begin{lem}\label{exchains}
Let $(M_k|k\in K)$ be a family of matroids, let $B_k\in \Ical(M_k)$ and let $C$ be a circuit for some $M_{k_0}$ such that $C \setminus B_{k_0}$ only contains one element, $e$.
If there is a $(B_k|k\in K)$-exchange chain from $x_0$ to $e$,
then for every $c\in C$, there is a $(B_k|k\in K)$-exchange chain from $x_0$ to $c$.
\end{lem}

\begin{proof}
Let $(y_0=x_0, k_0;y_1,k_1;\ldots;y_n=e)$ be an exchange chain from $x_0$ to $e$.
Then $(y_0=x_0, k_0;y_1,k_1;\ldots;y_n=e,k_0;c)$ is the desired exchange chain.
\end{proof}

\section{The Packing/Covering conjecture}

The matroid union theorem is a basic result in the theory of finite matroids. It gives a way to produce a new matroid $M = \bigvee_{k \in K}M_k$ from a finite family $(M_k | k \in K)$ of finite matroids on the same ground set $E$. We take a subset $I$ of $E$ to be $M$-independent iff it is a union $\bigcup_{k \in K}I_k$ with each $I_k$ independent in the corresponding matroid $M_k$. The fact that this gives a matroid is interesting, but a great deal of the power of the theorem comes from the fact that it gives an explicit formula for the ranks of sets in this matroid: 
\begin{equation}\label{rf0}
r_M(X) = \min_{X = P \dot{\cup} C} \sum_{k\in K} r_{M_k} (P)+|C|
\end{equation}
Here the minimisation is over those pairs $(P, C)$ of subsets of $X$ which partition $X$.

For infinite matroids, or infinite families of matroids, this theorem is no longer true \cite{union1}, in that $M$ is no longer a matroid. However, it turns out, as we shall now show, that we may conjecture a natural extension of the rank formula to infinite families of infinite matroids.

First, we state the formula in a way which does not rely on the assumption that $M$ is a matroid:
\begin{equation}\label{rf1}
\max_{ I_k\in \Ical(M_k)}\left|\bigcup_{k\in K} I_k\right| = \min_{E=P\dot\cup C}  \sum_{k\in K} r_{M_k} (P)+|C|
\end{equation}
Note that this is really only the special case of (\ref{rf0}) with $X = E$. However, it is easy to deduce the more general version by applying (\ref{rf1}) to the family $(M_k \restric_X | k \in K)$.

Note also that every element of the family over which the maximisation on the left is taken is at most as big as each member of the family over which the minimisation on the right is taken. To see this, note that $|\bigcup_{k\in K} (I_k\cap P)|\leq \sum_{k\in K} r_{M_k}(P)$ and $\bigcup_{k\in K} (I_k\cap C)\subseteq C$. So the formula is equivalent to the statement that we can find $(I_k | k \in K)$ and $P$ and $C$ with $P \dot \cup C = E$ so that
\begin{equation}\label{rf2}
\left|\bigcup_{k\in K} I_k\right| = \sum_{k\in K} r_{M_k} (P)+|C|\, .
\end{equation}
For this, what we need is to have equality in the two inequalities above, so we get 
\begin{equation}\label{rf3}
\left|\bigcup_{k\in K} (I_k\cap P)\right| = \sum_{k\in K} r_{M_k}(P)\text{ and }\bigcup_{k\in K} (I_k\cap C) = C \, .
\end{equation}
The equation on the left can be broken down a bit further: it states that each $I_k \cap P$ is spanning (and so a base) in the appropriate matroid $M_k \restric_P$, and that all these sets are disjoint. This is the familiar notion of a packing:

\begin{dfn}\label{pdef}
Let $(M_k | k \in K)$ be a family of matroids on the same ground set $E$. A  {\em packing} for this family consists of a spanning set $S_k$ for each $M_k$ such that the $S_k$ are all disjoint.
\end{dfn}

So the $I_k \cap P$ form a packing for the family $(M_k \restric_P | k \in K)$. In fact, in this case, each $I_k \cap P$ is a base in the corresponding matroid. In \autoref{pdef}, we do not require the $S_k$ to be bases, but of course if we have a packing we can take a base for each $S_k$ and so obtain a packing employing only bases.

Dually, the right hand equation in (\ref{rf3}) corresponds to the presence of a covering of $C$:

\begin{dfn}\label{cdef}
Let $(M_k | k \in K)$ be a family of matroids on the same ground set $E$. A  {\em covering} for this family consists of an independent set $I_k$ for each $M_k$ such that the $I_k$ cover $E$.
\end{dfn}

It is immediate that the sets $I_k \cap C$ form a covering for the family $(M_k \restric_C | k \in K)$. In fact we get the stronger statement that they form a covering for the family $(M_k.C | k \in K)$ where we contract instead of restricting, since for each $k$ we have that $I_k \cap P$ is an $M_k$-base for $P$, and we also have that $I_k$, which is the union of $I_k \cap C$ with $I_k \cap P$, is $M_k$-independent.

Putting all of this together, we get the following self-dual notion:

\begin{dfn}\label{pcdef}
Let $(M_k | k \in K)$ be a family of matroids on the same ground set $E$. We say this family {\em satisfies Packing/Covering} 
iff there is a partition of $E$ into two parts $P$ (called the {\em packing side}) and $C$ (called the {\em covering side}) such that 
$(M_k\restric_P | k \in K)$ has a packing, and $(M_k.C|k \in K)$ has a covering.
\end{dfn}

We have established above that this property follows from the rank formula for union, but the argument can easily be reversed to show that in fact Packing/Covering is equivalent to the rank formula, where that formula makes sense. However, Packing/Covering also makes sense for infinite matroids, where the rank formula is no longer useful. We are therefore led to the following conjecture:

\newtheorem*{thm3}{\autoref{pc}}

\begin{thm3}%\label{pc}
Every family of matroids on the same ground set satisfies Packing/Covering.
\end{thm3}

Because of this link to the rank formula, we immediately get a special case of this conjecture:

\begin{thm}
Every finite family of finite matroids on the same ground set satisfies Packing/Covering.
\end{thm}

Packing/Covering for pairs of matroids is closely related to another property which is conjectured to hold for all pairs of matroids.

\begin{dfn}
A pair $(M, N)$ of matroids on the same ground set $E$ {\em satisfies intersection} iff there is a subset $J$ of $E$, independent in both matroids, and a partition of $J$ into two parts $J^M$ and $J^N$ such that $$\Cl_M(J^M) \cup \Cl_N(J^N) = E\, .$$ 
\end{dfn}

\newtheorem*{thm2}{\autoref{int}}

\begin{thm2}%\label{int}
Every pair of matroids on the same ground set satisfies intersection.
\end{thm2}

We begin by demonstrating a link between Packing/Covering for pairs of matroids and intersection.

\begin{prop}\label{intispc}
Let $M$ and $N$ be matroids on the same ground set $E$. Then $M$ and $N$ satisfy intersection iff $M$ and $N^*$ satisfy Packing/Covering.
\end{prop}
\begin{proof}
Suppose first of all that $M$ and $N^*$ satisfy Packing/Covering, with packing side $P$ decomposed as $S^M \dot\cup S^{N^*}$ and covering side $C$ decomposed as $I^M \dot\cup I^{N^*}$. Let $J^M$ be an $M$-base of $S^M$, and $J^N$ an $N$-base of $C \setminus I^{N^*}$. $J = J^M \cup J^N$ is independent in $M$ since $J^N \subseteq I^M$ is independent in $M.C$ and $J^M$ is independent in $M \restric_P$. Similarly $J$ is independent in $N$ since $J^M \subseteq P \setminus S^{N^*}$ is independent in $N.P$ and $J^N$ is independent in $N \restric_C$. But also $$\Cl_M(J^M) \cup \Cl_N(J^N) = \Cl_M(S^M) \cup \Cl_N(C \setminus I^{N^*}) \supseteq P \cup C = E\, .$$

Now suppose instead that $M$ and $N$ satisfy intersection, as witnessed by $J = J^M \dot\cup J^N$. Let $J^M \subseteq P \subseteq \Cl_M(J^M)$ and $J^N \subseteq C \subseteq \Cl_N(J^N)$ be a partition of $E$ (this is possible since $\Cl_M(J^M) \cup \Cl_N(J^N) = E$). We shall show first of all that $M \restric_{P}$ and $N^* \restric_{P}$ have a packing, with the spanning sets given by $S^M = J^M$ and $S^{N^*} = P \setminus J^M$. $J^M$ is spanning in $M \restric_P$ since $P \subseteq \Cl_M(J^M)$, so it is enough to check that $P \setminus J^M$ is spanning in $N^* \restric_{P}$, or equivalently that $J^M$ is independent in $N.P$. But this is true since $J^N$ is an $N$-base of $C$ and $J^M \cup J^N$ is $N$-independent.

Similarly, $J^N$ is independent in $M.C$, and since $C \subseteq \Cl_N(J^N)$ $J^N$ is spanning in $N \restric_C$ and so $C \setminus J^N$ is independent in $N^*.C$. Thus the sets $I^M = J^N$ and $I^{N^*} = C \setminus J^N$ form a covering for $(M.C, N^*.C)$.
\end{proof}

\begin{cor}
If $M$ and $N$ are matroids on the same ground set then $M$ and $N$ satisfy Packing/Covering iff $M^*$ and $N^*$ do. \easy
\end{cor}

This corollary is not too hard to see directly. However, the following similar corollary is less trivial.

\newtheorem*{thm4}{\autoref{int_dual}}

\begin{thm4}
If $M$ and $N$ are matroids on the same ground set then $M$ and $N$ satisfy intersection iff $M^*$ and $N^*$ do. \easy
\end{thm4}

\autoref{intispc} shows that \autoref{int} follows from \autoref{pc}, but so far we can only use it to deduce \autoref{pc} for pairs of matroids from \autoref{int}. However, this turns out to be enough to give the whole of \autoref{pc}.

\begin{prop}\label{pc2ispc}
Let $(M_k|k \in K)$ be a family of matroids on the same ground set $E$, and let $M = \bigoplus_{k \in K}M_k$, on the ground set $E \times K$. Let $N$ be the matroid on the same ground set given by $\bigoplus_{e \in E} U_{1,K}^*$. Then the $M_k$ satisfy Packing/Covering iff $M$ and $N$ do.
\end{prop}
\begin{proof}
First of all, suppose that the $M_k$ satisfy Packing/Covering and let $P$, $C$, $S_k$ and $I_k$ be as in \autoref{pcdef}. We can partition $E \times K$ into $P' = P \times K$ and $C' = C \times K$. Let $S^M = \bigcup_{k \in K} S_k \times \{k\}$, and let $S^N = P' \setminus S^M$. $S^M$ is spanning in $M \restric_{P'}$ by definition, and since the sets $S_k$ are disjoint, there is for each $e \in P$ at most one $k \in K$ with $(e,k) \not \in S^N$. Thus $S^N$ is spanning in $N \restric_{P'}$. Similarly, let $I^M = \bigcup_{k \in K} I_k \times \{k\}$ and let $I^N = C' \setminus I^M$. $I^M$ is independent in $M.C'$ by definition, and since the sets $I_k$ cover $C$ there is for each $e \in E$ at least one $k \in K$ with $(e, k) \not \in I^N$. Thus $I^N$ is independent in $N.C'$.

Now suppose instead that $M$ and $N$ satisfy Packing/Covering, with packing side $P$ decomposed as $S^M \dot\cup S^{N}$ and covering side $C$ decomposed as $I^M \dot\cup I^{N}$. First we modify these sets a little so that the packing and covering sides are given by $\overline{P} \times K$ and $\overline{C} \times K$ for some sets $\overline{P}$ and $\overline{C}$. To this end, we let $\overline{P} = \{e \in E| (\forall k \in K) (e, k) \in P\}$, and $\overline{C} = \{e \in E | (\exists k \in K) (e, k) \in C\}$, so that $\overline{P}$ and $\overline{C}$ form a partition of $E$. Let $\overline{S}^N = S^N \cap (\overline{P} \times K)$ and $\overline{I}^N = I^N \cup ((\overline{C} \times K) \setminus C)$. We shall show that $(S^M, \overline{S}^N)$ is a packing for $(M\restric_{\overline{P} \times K}, N \restric_{\overline{P} \times K})$ and $(I^M, \overline{I}^N)$ is a covering for $(M.(\overline{C} \times K), N.(\overline{C} \times K))$.

 For any $e \in \overline{C}$, the restriction of the corresponding copy of $U_{1,K}^*$ to $P \cap (\{e\} \times K)$ is free, and so since the intersection of $S^N$ with this set is spanning there, it must contain the whole of $P \cap (\{e\} \times K)$. So since $S^M \subseteq P$ is disjoint from $S^N$, it can't contain any $(e, k)$ with $e \in \overline{C}$. That is, $S^M \subseteq \overline{P} \times K$. It also spans $\overline{P} \times K$ in $M$, since it spans the larger set $P$. For each $e \in \overline{P}$, $\overline{S}^N \cap (\{e\} \times K) = S^N \cap (\{e\} \times K)$ $N$-spans $\{e\} \times K$. Thus $\overline{S}^N$ $N$-spans $\overline{P} \times K$, so $(S^M, \overline{S}^N)$ is a packing for $(M\restric_{\overline{P} \times K}, N \restric_{\overline{P} \times K})$.

To show that $(I^M, \overline{I}^N)$ is a covering for $(M.(\overline{C} \times K), N.(\overline{C} \times K))$, it suffices to show that $\overline{I}^N$ is $N.(\overline{C} \times K)$-independent. For each $e \in \overline{C}$, the set $C \cap (\{e\} \times K)$ is nonempty, so the contraction of the corresponding copy of $U_{1,K}^*$ to this set consists of a single circuit, so there is some point in this set but not in $I^N$. Then that same point is also not in $\overline{I}^N$, and so $\overline{I}^N \cap 
(\{e\} \times K)$ is independent in the corresponding copy of $U_{1,K}^*$, so $\overline{I}^N$ is indeed $N.(\overline{C} \times P)$-independent.

Now that we have shown that $\overline{P} \times K$, $\overline{C} \times K$, $(S^M, \overline{S}^N)$ and $(I^M, \overline{I}^N)$ also witness that $M$ and $N$ satisfy Packing/Covering, we show how we can construct a packing and a covering for $(M_k \restric_{\overline{P}} | k \in K)$ and $(M_k . \overline{C} | k \in K)$ respectively.

For each $k \in K$ let $I_k = \{e \in E | (e, k) \in I^M\}$. Since, as we saw above, $I^M$ meets each of the sets $\{e\} \times K$ with $e \in \overline{C}$, the union of the $I_k$ is $\overline{C}$. Since also each $I_k$ is independent in $M_k.\overline{C}$, they form a covering for $(M_k.\overline{C} | k \in K)$. Similarly, let $S_k = \{e \in E | (e,k) \in S^M\}$. Since the intersection of $\overline{S}^N$ with $\{e\} \times K$ is spanning in the corresponding copy of $U_{1,k}^*$ for any $e \in \overline{P}$, it follows that for such $e$ it misses at most one point of this set, so that there can be at most one point in $S^M \cap (\{e\} \times K)$, so the $S_k$ are disjoint. Thus they form a packing of $(M_k\restric_{\overline{P}} | k \in K)$.
\end{proof}

\begin{cor}
The following are equivalent:
\begin{enumerate}
\item Intersection holds for any pair of matroids (\autoref{int}).
\item Intersection holds for any pair of matroids in which the second is a direct sum of copies of $U_{1,2}$.
\item Packing/Covering holds for any pair of matroids.
\item Packing/Covering holds for any pair of matroids in which the second is a direct sum of copies of $U_{1,2}$.
\item Packing/Covering holds for any family of matroids (\autoref{pc}).
\end{enumerate}
\end{cor}
\begin{proof}
We shall prove the following equivalences.
$$\xymatrix{2 \ar@{<->}[r] & 4 \ar@{<->}[d] &\\ 1 \ar@{<->}[r] & 3 \ar@{<->}[r] & 5}$$
The equivalences of (1) with (3) and (2) with (4) both follow from \autoref{intispc}. (3) evidently implies (4), but we can also get (4) from (3) by applying \autoref{pc2ispc}. Similarly, (5) evidently implies (3) and we can get (5) from (3) by applying \autoref{pc2ispc}.
\end{proof}

\section{A special case of the Packing/Covering conjecture}

In \cite{Aharoni:Ziv:92}, Aharoni and Ziv prove a special case of the intersection conjecture. Here we employ a simplified form of their argument to prove a special case of the Packing/Covering conjecture. Our simplification also yields a slight strengthening of their theorem. 

Key to the argument is the notion of a wave.

\begin{dfn}
Let $(M_k|k \in K)$ be a family of matroids all on the ground set $E$. A {\em wave} for this family is a subset $P$ of $E$ 
together with a packing $(S_k | k \in K)$ of $(M_k \restric_P|k \in K)$. In a slight abuse of notation, we shall sometimes refer to the wave just as $P$ or say that elements of $P$ are in the wave. 
A wave is a {\em hindrance} if the $S_k$ don't completely cover $P$. The family is {\em unhindered} if there is no hindrance, and {\em loose} if the only wave is the empty wave.
\end{dfn}

\begin{comment}COMMENT TO DEF OF HINDRANCE:
Ellie Berger said that a slightly different notation of hindrance might be the right one:
A family $S_k$ is a hindrance if $S_k\in cl_{M_l} S_l$ for all $k,l$ and there is $e\notin \bigcup S_k$ with $e\in cl_{M_k} S_k$ for some $k$.
With this notation there are more hindrances (any family with a loop in any $
M_k$ is already hindered). However, for finite matroids our notation is strong enough to prove pack-cov. 
\end{comment}

%With a slight abuse of notation, we say that a set $P$ is a wave if there is a family of disjoint spanning sets $(S_k | k \in K)$ of $(M_k \restric_P|k \in K)$ .

\begin{rem}
Those familiar with Aharoni and Ziv's notion of wave should observe that if $(P, (S_1, S_2))$ is a wave as above and we let $F$ be an $M_2$-base of $S_2$ then $F$ is not only $M_2$-independent but also $M_1^*.P$-independent, since $S_1 \subseteq P \setminus F$ is $M_1 \restric_P$-spanning. Now since $P \subseteq \Cl_{M_2}(F)$, we get that $F$ is also $M_1^*.\Cl_{M_2}(F)$-independent. Thus $F$ is a wave in the sense of Aharoni and Ziv for the matroids $M_1^*$ and $M_2$. There is a similar correspondence of the other notions defined above. 

Similarly, they say that the pair $(M_1, M_2)$ is {\em matchable} iff there is a set which is $M_1$-spanning and $M_2$-independent. Those interested in translating between the two contexts should note that there is a covering for $(M_1, M_2)$ iff $(M_1^*, M_2)$ is matchable.
\end{rem}

We define a partial order on waves by $(P, (S_k | k \in K)) \leq (P', (S_k'  | k \in K))$ iff $P \subseteq P'$ and for each $k \in K$ we have $S_k \subseteq S_k'$. We say a wave is {\em maximal} iff it is maximal with respect to this partial order. 

\begin{lem}\label{joinwaves}
Let $(M_k | k \in K)$ be a family of matroids on the same ground set $E$, and let $((P^{\beta}, (S^{\beta}_k | k \in K))|\beta < \alpha)$ a family of waves indexed by some ordinal $\alpha$. Then there is a wave $(P, (S_k, | k \in K))$ with $P = \bigcup_{\beta < \alpha} P^{\beta}$ and $P \geq P_0$. 
\end{lem}
\begin{proof}
For each $\beta < \alpha$, let $Y^{\beta} = P^{\beta} \setminus \bigcup_{\gamma < \beta} P^{\gamma}$. For $k \in K$, let $S_k = \bigcup_{\beta < \alpha}(Y^{\beta} \cap S^{\beta}_k)$. These are clearly disjoint subsets of $P$: we aim to show that they form a packing. We shall show by induction on $\beta < \alpha$ that for each $k \in K$ we have $P^{\beta} \subseteq \Cl_{M_k}(S_k)$. By the induction hypothesis, we have that $S_k^{\beta} \setminus Y^{\beta} \subseteq \bigcup_{\gamma < \beta} P^{\gamma} \subseteq \Cl_{M_k}(S_k)$, so $P^{\beta} \subseteq \Cl_{M_k}(S^{\beta}_k) \subseteq \Cl_{M_k}(\Cl_{M_k}(S_k)) = \Cl_{M_k}(S_k)$. 

It follows that $P \subseteq \Cl_{M_k}(S_k)$, so the $S_k$ form a packing for $(M_k \restric_P)$ as desired.
\end{proof}

\begin{cor}\label{maxwave}
For any wave $P$ there is a maximal wave $P_{\max} \geq P$. 
\end{cor}
\begin{proof}
We apply \autoref{joinwaves} to a family of waves with $P$ as the first element and which includes all waves.
\end{proof}

\begin{cor}\label{coverall}
If $P_{\max}$ is a maximal wave then anything in any wave $P$ is in $P_{\max}$.
\end{cor}
\begin{proof}
We apply \autoref{joinwaves} to the pair $(P_{\max}, P)$.
\end{proof}

\begin{lem}\label{span_edge}
For any $e\in E$, any maximal wave $P$ satisfies $e\in cl_{M_k} P$ whenever there is any wave $P'$ with $e\in cl_{M_k} P'$.

In particular, if $e$ is not contained in any wave, there are at least two $k$ such that, for every wave $P'$,  $e\notin cl_{M_k} P'$.
\end{lem}

%\begin{lem}\label{span_edge}
%Let $e\in E$. There is a wave $P$ such that for every other wave $P'$ and every $k$ with $e\in cl_{M_k} P'$ implies
%$e\in cl_{M_k} P$.
%In particular, if $e$ is not contained in any wave, there are at least two $k$ such that $e\notin cl_{M_k} P'$ for every wave $P'$.
%\end{lem}

\begin{proof}
Let $(P,(S_k|k\in K)) $ be a maximal wave.
By \autoref{coverall} for any wave $(P',(S_k'|k\in K))$ we have $S_k'\subseteq cl_{M_k}S_k$.
Thus $e\in cl_{M_k} P'=cl_{M_k}S_k'$ implies $e\in cl_{M_k} P$, as desired.

For the second assertion, assume toward contradiction that there is at most one $k_0$ such that, for every wave $P'$,  $e\notin cl_{M_{k_0}} P'$.
Then $e\in cl_{M_{k}} P$ for all $k\neq k_0$.
But then the following is a wave and contains $e$: \\
$X:=(P+e,(\overline{S}_k|k \in K))$ where $\overline{S}_{k_0} = S_{k_0} + e$ and $\overline{S}_k = S_k$ for other values of $k$.
 This is a contradiction.
\end{proof}

\begin{lem}
Let $(P, (S_k | k \in K))$ be a wave for a family $(M_k| k \in K)$ of matroids. 
Let $(P', (S'_k | k \in K))$ be a wave for the family $(M_k/P | k \in K)$. 
Then $(P \cup P', (S_k \cup S_k' | k \in K))$ is a wave for the family $(M_k| k \in K)$. If either $P$ or $P'$ is a hindrance then so is $P \cup P'$.
\end{lem}
\begin{rem}
In fact, though we will not need this, a similar statement can be shown for an ordinal indexed family of waves $P^{\beta}$, with $P^{\beta}$ a wave for the family $(M_k/\bigcup_{\gamma < \beta} P^{\gamma}| k \in K)$.
\end{rem}
\begin{proof}
By \autoref{sthg}, each $S_k \cup S_k'$ spans $P \cup P'$, and they are clearly disjoint. 
If the $S_k$ don't cover some point of $P$ then the $S_k \cup S_k'$ also don't cover that point, and the argument in the case where $P'$ is a hindrance is similar.
\end{proof}

\begin{cor}\label{loose}
For $P_{\max}$ as in \autoref{maxwave}, the family $(M_k / P_{\max} | k \in K)$ is loose.
\end{cor}

We are now in a position to present another Conjecture equivalent to the Packing/Covering Conjecture. It is for this new form that we shall present our partial proof.

\begin{con}\label{unhcov}
Any unhindered family of matroids has a covering.
\end{con}

\begin{prop}\label{unhindcov_is_pc}
\autoref{unhcov} and \autoref{pc} are equivalent.
\end{prop}
\begin{proof}
First of all, suppose that \autoref{pc} holds, and that we have an unhindered family $(M_k | k \in K)$ of matroids. Using \autoref{pc}, we get $P$, $C$, $S_k$ and $I_k$ as in \autoref{pcdef}. Then $(P, (S_k | k \in K))$ is a wave, and since it can't be a hindrance the sets $S_k$ cover $P$. They must also all be independent, since otherwise we could remove a point from one of them to obtain a hindrance. So the sets $S_k \cup I_k$ give a covering for $(M_k | k \in K)$. 

Now suppose instead that \autoref{unhcov} holds, and let $(M_k | k \in K)$ be any family of matroids on the ground set $E$. Then let $(P, (S_k | k \in K))$ be a maximal wave, as in \autoref{maxwave}. By \autoref{loose}, $(M_k/P | k \in K)$ is loose, and so in particular this family is unhindered. So it has a covering $(I_k | k \in K)$. Taking covering side $C = E \setminus P$, this means that the $M_k$ satisfy Packing/Covering.
\end{proof}

\begin{lem}\label{edge_simply}
Suppose that we have an unhindered family $(M_k|k \in K)$ of matroids on a ground set $E$. Let 
$e \in E$ and $k_0 \in K$ such that for every wave $P$ we have $e\notin cl_{M_{k_0}} P$. Then the family $(M_k' | k \in K)$ on 
the ground set $E - e$ is also unhindered, where $M_{k_0}' = M_{k_0} / e$ but $M_k' = M_k \backslash e$ for other values of $k$.
\end{lem}

\begin{proof}

Suppose not, for a contradiction, and let $(P, (S_k | k \in K))$ be a hindrance for 
$(M_k' | k \in K)$. 
Without loss of generality, we assume that the $S_k$ are bases of $P$.
Let $\overline{S}_k$ be given 
by $\overline{S}_{k_0} = S_{k_0} + e$ and $\overline{S}_k = S_k$ for other values of $k$. 
Note that $\overline{S}_{k_0}$ is independent because otherwise $e\in cl_{M_{k_0}}P$.
Let $P'$ be the set of $x \in P$ such that there is no $(\overline{S}_k| k \in K)$-exchange chain from $x$ to $e$.

Let $x_0 \in P \setminus \bigcup_{k \in K} S_k$.  If $x_0 \in P'$, then we will show that 
$(P', P' \cap \overline{S}_k)$ is a wave containing $x_0$. 
This contradicts the assumption that $(M_k | k \in K)$ is unhindered. Since $e\notin P'$, we have $P' \cap\overline{S}_k=P' \cap S_k$ for every $k$.
So it suffices to show for every $k$ that every $x\in P'\setminus P' \cap \overline{S}_k$ is $M_k$-spanned by $P' \cap \overline{S}_k$. 
Let $C$ be the unique circuit contained in $x+\overline{S}_k$.
If $x\in P'$, then $C\subseteq P'$ by \autoref{exchains}, so $x\in  cl_{M_{k}}P'\cap \overline{S}_k$, as desired.

If $x_0\notin P'$, there is a shortest $(\overline{S}_k|k\in K)$-exchange chain 
$(y_0=x_0, k_0;y_1,k_1;\ldots;y_n=e)$ from $x_0$ to $e$.
Let $\overline{S}_k':=\overline{S}_k\cup \{y_l|k_{l}=k\} \setminus \{y_{l+1}|k_{l}=k\}$.
By \autoref{runchains}, $\overline{S}_k'$ is $M_k$-independent and $cl_{M_{k}}\overline{S_{k}}=cl_{M_{k}}\overline{S_{k}}'$ for all $k\in K$.
Thus each $\overline{S}_k$ $M_k$-spans $P$ but avoids $e$, in other words:  $(P,(\overline{S}_k'|k\in K))$ is an $(M_k|k\in K)$-wave.
But also $e\in  cl_{M_{k_0}}P$ since $e\in \overline{S}_{k_0}$, a contradiction.
\end{proof}

We will now discuss those partial versions of \autoref{unhcov} which we can prove. We would like to produce a covering of the ground set by independent sets - and that means that we don't want any of the sets in the covering to include any circuits for the corresponding matroid. First of all, we show that we can at least avoid {\em some} circuits. In fact, we'll prove a slightly stronger theorem here, showing that we can specify a countable family of sets, which are to be avoided whenever they are dependent. In all our applications, the dependent sets we care about will be circuits.

\begin{thm}\label{cutcount}
Let $(M_k|k \in K)$ be an unhindered family of matroids on the same ground set $E$. 
Suppose that we have a sequence of subsets $o_n$ of $E$. Then there is a family $(I_k| k \in K)$ covering $E$ such that for no $k \in K$ and $n \in \Nbb$ do we have both $o_n \subseteq I_{k}$ and $o_n$ dependent in $M_k$.
\end{thm}

%TODO made stronger, proof changed

\begin{proof}
If some wave includes the whole ground set, then as the family is unhindered, this wave would yield the desired covering.
Unfortunately, we may not assume this. Instead, we recursively build a family $(J_k|k\in K)$ of disjoint sets such that some wave $(P, (S_k|k\in K))$
for the  $M_k/J_k-\bigcup_{l\neq k} J_l$ includes enough of $E-\bigcup_{k}J_k$ that any family $(I_k|k \in K)$ covering $E$ and with 
$I_k \cap (P \cup \bigcup_{k \in K} J_k) = S_k\cup J_k$ will work.

We construct $J_k$ as the nested union of some $(J_k^n|n\in \Nbb)$ with the following properties.
Abbreviate $M_k^n:=M_k/J_k^n-\bigcup_{l\neq k} J_l^n$.

\begin{enumerate}
\item $J_k^n$ is independent in $M_k$. \label{prf1}
\item For different $k$, the sets $J_k^n$ are disjoint. \label{prf1.5}
\item $(M_k^n|k\in K)$ is unhindered. \label{prf2}
\item Either the set $o_n-\bigcup_{k\in K} J_k^n$ is included in some $(M_k^n|k\in K)$-wave
or there are distinct $l, l'$ such that there is some $e\in o_n\cap J_{l}^n$ and some $e'\in o_n\cap J_{l'}^n$.\label{prf3}
%or if $M_k$ witnesses $o_i$, then there is some $e\in o_{i}$ such that either $e\cup J_k^{i-1}$ is dependent or $e\in J_{k'}^n$ with $k'\neq k$. \label{prf3}
\end{enumerate}

%If $(o_i| i\in \Nbb)$ is finite, then we will get stationary after the last $o_i$.
Put $J_k^0:=\emptyset$ for all $k$. Assume that we have already constructed $J_k^n$ satisfying (\ref{prf1})-(\ref{prf3}).

If (\ref{prf3}) with $o_{n+1}$ in place of $o_n$ is already satisfied by the  $(J_k^n|k\in K)$ we can simply take $J_k^{n+1}:=J_k^n$ for all $k$.

Otherwise by \autoref{maxwave}, there is some $e\in o_{n+1}-\bigcup_{k\in K} J_k^n$ not in any $(M_k^n|k\in K)$-wave.
%Then we further have that $e\cup J_w^{n}$ is independent.
By \autoref{span_edge}, there are at least two $k \in K$ such that $e\notin cl_{M_k} P'$ for every wave $P'$.
In particular, $e$ is not a loop ($\{e\}$ is independent) in $M_k$ for those two $k$.
Let $l$ be one of these two values of $k$.
Now let $\overline{J_l^{n+1}}:=J_l^n+e$ and $\overline{J_k^{n+1}}:=J_k^n$ for $k\neq l$.
Then the $\overline{J_l^{n+1}}$ satisfy (\ref{prf1})-(\ref{prf1.5}).
By \autoref{edge_simply} and the choice of $e$, we also have (\ref{prf2}).

If the $\overline{J_l^{n+1}}$ already satisfy (\ref{prf3}), then the are done.
Else, to obtain (\ref{prf3}), repeat the induction step so far and find $e'\in o_{n+1}-\bigcup_{k\in K} \overline{J_k^{n+1}}-e$ not in any $(\overline{M_k^n}|k\in K)$-wave.
Here $\overline{M_k^n}$ is $M_k^n/e$ if $k= l$ and $M_k^n-e$ otherwise.
Further we find, $l'\neq l$ such that $e'$ is independent in $\overline{M_{l'}^n}$ and $e'\notin cl_{M_l} P'$ for every wave $P'$.
Now let $J_{l'}^{n+1}:=\overline{J_{l'}^{n+1}}+e'$ and $J_k^{n+1}:=\overline{J_k^{n+1}}$ for $k\neq l'$.
Then the ${J_k^{n+1}}$ satisfy (\ref{prf1})-(\ref{prf1.5}) and now also (\ref{prf3}).
By \autoref{edge_simply} and the choice of $e'$, we also have (\ref{prf2}).

We now define a new family of matroids by $M_k' := M_k/J_k-\bigcup_{l\neq k} J_l$, and we construct an $(M_k'|k\in K)$-wave $(P, (S_k|k\in K))$. 
We once more do this by taking the union of a recursively constructed nested family. Explicitly, we take $S_k = \bigcup_{n \in \Nbb} S_k^n$ and 
$P = \bigcup_{n \in \Nbb} P^n$, where for each $n$ the wave $W_n = (P^n, (S_k^n| k \in K))$ is a maximal wave for $(M_k^n | k \in K)$ and the $S_k^n$ are nested. 
We can find such waves using \autoref{maxwave}: for each $n$ we have that $W^n$ is also a wave for $(M^{n+1}_k | k \in K)$ since in our construction 
we never contract or delete anything which is in a wave.
%Since in the construction of the $J_k$ we only contracted elements not contained in a wave, every $(M_k^n|k\in K)$-wave is still an $(M_k/J_k-\bigcup_{l\neq k} J_l|k\in K)$-wave.

Now let $(I_k | k \in K)$ be chosen so that  $\bigcup I_k=E$ and for each $k$ we have $I_k \cap (P \cup \bigcup_{k \in K} J_k) = S_k\cup J_k$. 
Suppose for a contradiction that for some pair $(k_0,n)$ we have $o_n \subseteq I_{k_0}$ and $o_n$ is dependent in $M_{k_0}$. 
Then by (\ref{prf3}), either the set $o_n-\bigcup_{k\in K} J_k^n$ is included in some $(M_k^n|k\in K)$-wave
or there are distinct $l, l'$ such that there is some $e\in o_n\cap J_{l}^n$ and some $e'\in o_n\cap J_{l'}^n$. 
In the second case, clearly $o_n \nsubseteq I_{k_0}$.

In the first case, we will find a hindrance for $(M_k^n|k\in K)$, which contradicts (\ref{prf2}).
It suffices to show that $S_{k_0}^n$ is dependent.
As $o_n\subseteq I_{k_0}$, we have $o_n-J_{K_0}^n\subseteq S_{k_0}^n$.
Note that $o_n-J_{K_0}^n$ is non-empty by (\ref{prf1}).
But now $o_n-J_{K_0}^n$ is dependent in $M_{k_0}^n$ and contained in $S_{k_0}^n$, a contradiction.

\begin{comment}
 In the first case, the set $o_n-\bigcup_{k\in K} J_k^n$ can't be empty by (\ref{prf1}). 
It can't meet any $S_k^n$ with $k \neq k_0$, so either it is a nonempty subset of $S^n_{k_0}$ or else the $S_k^n$ don't cover $P^n$. 
In either case, we get a hindrance for the $M^n_k$ (if it is a nonempty subset of $S_{k_0}^n$ then, since it is $M_{k_0}^n$-dependent we can remove 
some element from $S^n_{k_0}$ and still have a spanning set). This contradicts (\ref{prf2}).
\end{comment}

\end{proof}

Note that, in particular, if we have a countable family of matroids each with only countably many circuits then \autoref{cutcount} applies in order to prove \autoref{pc} in that special case. Requiring only countably many circuits might seem quite restrictive, but there are many cases where it holds:

\begin{prop}
A matroid of any of the following types on a countable ground set has only countably many circuits:
\begin{enumerate}
\item A finitary matroid.\label{cc1}
\item A matroid whose dual has finite rank.\label{cc2}
\item A direct sum of matroids each with only countably many circuits.\label{cc3}
\end{enumerate}
\end{prop}
\begin{proof}
(\ref{cc1}) follows from the fact that the countable ground set has only countably many finite subsets.
For (\ref{cc2}), since every base $B$ has finite complement, there are only countably many bases.
As every circuit is a fundamental circuit for some base, there can only be countably many circuits, as desired.
For (\ref{cc3}), there can only be countably many nontrivial summands in the direct sum since the ground set is countable, and the result follows.
\end{proof}

In particular, \autoref{cutcount} applies to any countable family of matroids each of which is a direct sum of matroids that are finitary or whose duals have finite rank. This includes the main result of Aharoni and Ziv in \cite{Aharoni:Ziv:92}, if the ground set if $E$ is countable, by \autoref{intispc}.

If we have a family of sets $(I_k | k \in K)$ which does not form a covering, because some elements aren't independent, how might we tweak it to make them more independent? Suppose that the reason why $I_k$ is dependent is that it contains a circuit $o$ of $M_k$, but that $o$ also includes a bond for another matroid $M_{k'}$ from our family. Then we could move some point from $I_k$ into $I_{k'}$ to remove this dependence without making $I_{k'}$ any more dependent%
\footnote{Note that wlog we may assume that the $I_k$ are disjoint. Then any new circuits in $I_{k'}$ would have to meet the bond in just one point, which is impossible by \autoref{cir_cocir}.}. We are therefore not so worried about circuits including bonds in this way as we are about other sorts of circuits. Therefore we now consider cases where most circuits do include such bonds: 

\begin{comment}
To find a covering for the family $(M, M^*)$ is easy: Any base of $M$ together with its complement form a covering. We might hope that if we had $M$ and $N$ with $N$ sufficiently close to $M^*$ then we could also find a covering. Our next result includes a condition like this: that each matroid in the family we are considering should be, in a sense, mostly covered by the duals of the other matroids.
\end{comment}

%\newcommand{\cw}{at most countably weird} already defined above

\begin{dfn}\label{weird}
Let $(M_k|k \in K)$ be a family of matroids on the same ground set $E$. For each $k \in K$ we let $W_k$ be the set of all $M_k$-circuits that do not contain an $M_{k'}$-bond with $k' \not = k$.
Call the family $(M_k|k \in K)$ of matroids \emph{\cw} if $\bigcup W_k$ is at most countable.%
\end{dfn}

Note that if $E$ is countable then $(M_k|k \in K)$ is \cw\ if and only if $\bigcup W_k^\infty$ is countable where $W_k^\infty$ is the subset of $W_k$ consisting only 
of the infinite circuits in $W_k$.

\begin{thm}\label{unhindcov}
Any unhindered and \cw\ family $(M_k|k \in K)$ of matroids has a covering.
\end{thm}

\begin{proof}

Apply \autoref{cutcount} to $(M_k|k \in K)$ where the $o_n$ enumerate $\bigcup W_k$ where the $W_k$ are defined as in \autoref{weird}.

So far $(I_k|k\in K)$ is not necessarily a covering since each $I_k$ might still contain circuits. But by the choice of the family of circuits each circuit contained in $I_k$ 
contains an $M_{k'}$-bond with $k' \not = k$. 

In the following, we tweak $(I_k|k\in K)$ to obtain a covering $(L_k|k\in K)$.
First extend $I_k$ into a minimal $M_k$-spanning set $B_k$ by $(IM)^*$.
We obtain $L_k$ from $B_k$ by removing all elements in $I_k \cap \bigcup_{l\neq k} B_l$.
We can suppose without loss of generality $(I_k|k\in K)$ was a partition of $E$, and so the  family $(L_k|k\in K)$ covers $E$. 
It remains to show that $L_k$ is independent.
For this, assume for a contradiction that $L_k$ contains an $M_k$-circuit $C$.
By the choice of $B_k$, the circuit $C$ is contained in $I_k$.
In particular, $C$ contains an $M_l$-bond $X$ for some $l\neq k$.
By construction $B_l$ meets $X$ and thus $C$.
As $C\subseteq I_k$, the circuit $C$ is not contained in $L_k$, a contradiction.
So $(L_k|k\in K)$ is the desired covering.
\end{proof}

%Old thm, I think that we cant proof this one
%\begin{thm}
%Suppose the family $(M_k|k \in K)$ of matroids on a ground set $E$ is unhindered and there are only countably many $o$ which are circuits in some $M_k$ but contain no bond of any $M_{k'}$ with $k' \not = k$. Then there is a covering for this family.
%\end{thm}

We can now apply the argument of \autoref{unhindcov_is_pc} to obtain the following:

\begin{cor}\label{pack_cov_thm}
Any \cw\ family $(M_k|k \in K)$ of matroids  satisfies Packing/Covering. \easy
%Def of cw changed on a countable ground set $E$
\end{cor}

However, there are still some important open questions here.

\begin{dfn}[\cite{union2}]
The \emph{finitarization of a matroid $M$} is the matroid $M^{fin}$
whose circuits are precisely the finite circuits of $M$%
\footnote{It is easy to check that $M^{fin}$ is indeed a matroid \cite{union2}}.
A matroid is called \emph{nearly finitary} if every base misses at most finitely elements of some base of the finitarization. 
\end{dfn}

From \autoref{intispc} and the corresponding case of matroid intersection \cite{union2} we obtain the following: 

\begin{cor}\label{int_nf}
 The Packing/Covering conjecture is true for two nearly finitary matroids. 
\end{cor}

By \autoref{pc2ispc} \autoref{int_nf} implies the packing covering conjecture for finite families of nearly finitary matroids.
We do not know the answer to the following question.

\begin{oque}
 Is the Packing/Covering Conjecture true for any (countably) infinite family of nearly finitary matroids?
\end{oque}

In a similar way, we have the following question.

\begin{oque}\label{o_fini}
 Is the Packing/Covering Conjecture true for arbitrary families of finitary matroids?
\end{oque}

%TODO: Can we get $\infty$ many nearly finitary?

\section{Base covering}\label{sec:cov}

The well-known base covering theorem reads as follows.

\begin{thm}\label{cov_thm}
Any family of finite matroids $(M_k| k\in K)$ on a finite common ground set $E$ has a covering if and only if for every finite set $X\subseteq E$ the following holds.
\[
\sum_{k\in K} r_{M_k}(X)\geq |X|
\]
\end{thm}

Taking the family to contain only one matroid, consisting of one infinite circuit, we see that this theorem does not extend verbatim to infinite matroids.
However, \autoref{cov_thm} extends verbatim to finite families of finitary matroids by compactness \cite{union1}%
\footnote{The argument in \cite{union1} is only made in the case where all $M_k$ are the same but it easily extends to finite families of arbitrary finitary matroids.}. 
The requirement that the family is finite is necessary as $(U_k=U_{1,\Rbb}|k\in \Nbb)$ satisfies the rank formula but does not have a covering.

In the following, we conjecture an extension of the finite base covering theorem to arbitrary infinite matroids.
Our approach is to replace the rank formula by a condition that for finite sets $X$ is implied by the rank formula but is still meaningful for infinite sets.
A first attempt might be the following:

\begin{equation}\label{cov_cond_pre}
\begin{minipage}[c]{0.8\textwidth}
Any packing for the family $(M_k\restric_X | k\in K)$ is already a covering.
\end{minipage}
\end{equation}

\begin{figure}[htbp]
	\centering
	\subfloat{\includegraphics[width=9cm]{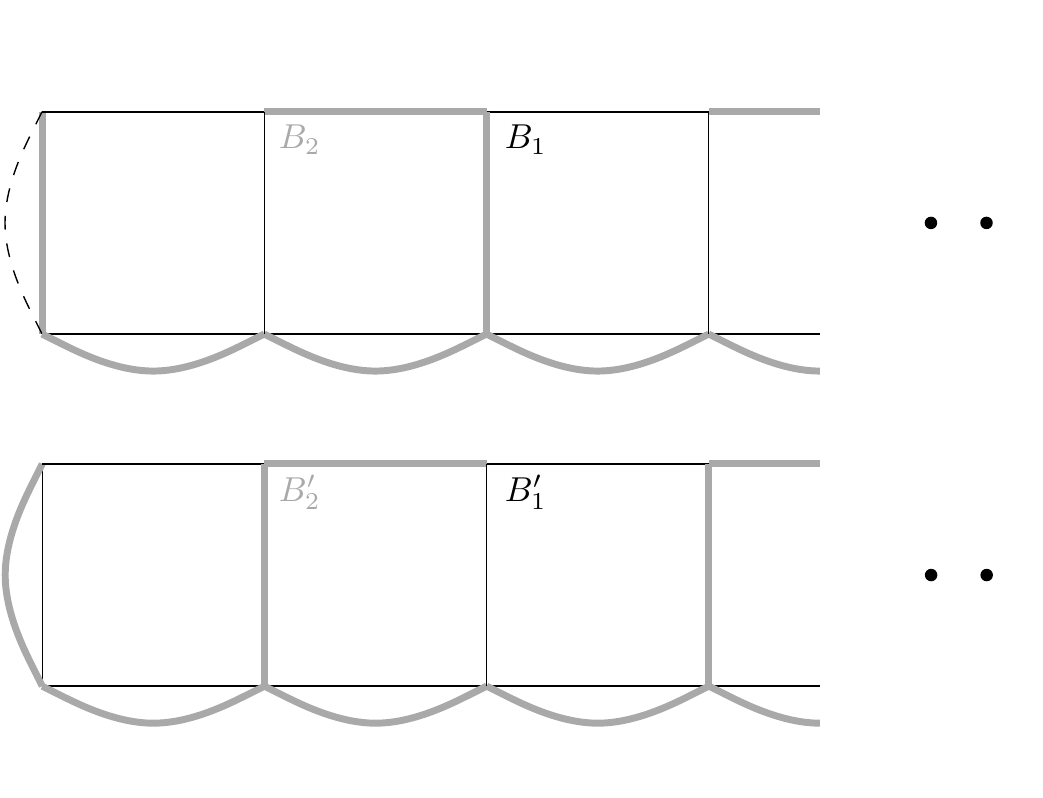}}
	\caption{Above is a base packing which isn't a base covering. Below that is a base covering for the same matroids, namely the finite cycle matroid for the graph, taken twice.}
	\label{fig:bases}
\end{figure}

Indeed, for finite $X$, if $(M_k\restric_X | k\in K)$ has a packing and there is an element of $X$ not covered by the corresponding bases, then this violates the rank formula.
However, there are infinite matroids that violate (\ref{cov_cond_pre}) and still have a covering, see \autoref{fig:bases}.

We propose to use instead the following weakening of (\ref{cov_cond_pre}).
\begin{equation}\label{cov_cond}
\begin{minipage}[c]{0.8\textwidth}
If $(M_k\restric_X | k\in K)$ has a packing,
then it also has a covering.
\end{minipage}
\end{equation}

To see that (\ref{cov_cond}) does not imply the rank formula for some finite $X$, consider the family $(M,M)$, where $M$ is the finite cycle matroid of the graph
$$\xymatrix{\bullet \ar@{-}[r] & \bullet \ar@{-}@/^5pt/[r] \ar@{-}@/_15pt/[r] \ar@{-}@/^15pt/[r] \ar@{-}@/_5pt/[r] & \bullet}$$ 
This graph has an edge not contained in any cycle (so that $(M,M)$ does not have a packing) but enough parallel edges to make the rank formula false.

Using (\ref{cov_cond}), we obtain the following:

\begin{con}[Covering Conjecture]\label{cov_conj}
A family of matroids $(M_k|k\in K)$ on the same ground set $E$ has a covering 
if and only if (\ref{cov_cond}) is true for every $X\subseteq E$.
\end{con}

\begin{prop}\label{pc_is_cov_conj}
\autoref{pc} and \autoref{cov_conj} are equivalent.
\end{prop}

\begin{proof}
For the ``only if'' direction, note that \autoref{cov_conj} implies \autoref{unhcov}, which by \autoref{unhindcov_is_pc} implies \autoref{pc}.

For the ``if'' direction, note that by assumption we have a partition $E=P\dot\cup C$ such that
there exists disjoint $M_k\restric_{P}$-spanning sets $S_k$ and $M_k.C$-independent sets $I_k$ covering $C$.
By (\ref{cov_cond}), $( M_k\restric_{P}|k\in K)$ has a covering with sets $B_k$, where $B_k\in \Ical(M_k\restric_{P})$.
As $I_k\cup B_k\in \Ical(M_k)$, the sets $I_k\cup B_k$ form the desired covering.
\end{proof}

As Packing/Covering is true for finite matroids, \autoref{pc_is_cov_conj} implies the non-trivial direction of \autoref{cov_thm}.
By \autoref{pack_cov_thm} we obtain the following applications.

\begin{cor}
A family of matroids $(M_k|k\in K)$ as in \autoref{pack_cov_thm} has a covering if and only if (\ref{cov_cond}) is true for every $X\subseteq E$.
\end{cor}

Let us now specialise to graphs.

\begin{dfn}
The bases of the topological cycle
matroid are called \emph{topological trees} and the bases of the algebraic cycle
matroid are called \emph{algebraic trees}.
Using this we define \emph{topological tree-packing, topological tree-covering, algebraic tree-packing, algebraic tree-covering}.
\end{dfn}

\begin{comment}
\begin{cor}[Base covering for arbitrary families of the finite circuit matroid]
A family of graphs $(G_k|k\in K)$ on a common ground set $E$ has a tree-covering if and only if the following is true for every $X\subseteq E$.
\begin{equation}\label{cov_cond_fc}
\begin{minipage}[c]{0.8\textwidth}
If $(G_k[X] | k\in K)$ has a tree-packing,
then it also has a tree-covering.
\end{minipage}
\end{equation}
\end{cor}
\end{comment}

\begin{cor}[Base covering for the topological cycle matroids]
A family of graphs $(G_k|k\in K)$ with a common edge set $E$ has a topological tree-covering if and only if the following is true for every $X\subseteq E$.
\begin{equation}\label{cov_cond_tc}
\begin{minipage}[c]{0.8\textwidth}
If $(G_k[X] | k\in K)$ has a topological tree-packing,
then it also has a topological tree-covering.
\end{minipage}
\end{equation}
\end{cor}

\begin{cor}[Base covering for the algebraic cycle matroids of locally finite graphs]
A family of locally finite graphs $(G_k|k\in K)$ with a common edge set $E$ has an algebraic tree-covering if and only if the following is true for every $X\subseteq E$.
\begin{equation}\label{cov_cond_ac}
\begin{minipage}[c]{0.8\textwidth}
If $(G_k[X] | k\in K)$ has an algebraic tree-packing,
then it also has an algebraic tree-covering.
\end{minipage}
\end{equation}
\end{cor}

\section{Base packing}

The well-known base packing theorem reads as follows.

\begin{thm}\label{pack_thm}
Any family of finite matroids $(M_k| k\in K)$ on a finite common ground set $E$ has a packing if and only if for every finite set $Y\subseteq E$ the following holds.
\[
\sum_{k\in K}  r_{M_k.Y}(Y) \leq |Y|
\]
\end{thm}

Aigner-Horey, Carmesin and Fr\"ohlich \cite{union1} extended this theorem to families consisting of finitely many copies of the same co-finitary matroid. 
We extend this to arbitrary co-finitary families.

\begin{thm}\label{pack_thm_fam}
Any family of co-finitary matroids $(M_k| k\in K)$ on a common ground set $E$ has a packing if and only if for every finite set $Y\subseteq E$ the following holds.
\[
\sum_{k\in K}  r_{M_k.Y}(Y) \leq |Y|
\]
\end{thm}

\begin{proof}[Proof by a compactness argument.]
We will think of partitions of the ground set $E$ as functions from $E$ to $K$ - such a function $f$ corresponds to a partition $(S^f_k | k \in K)$, given by $S_k^f= \{e \in E | f(e) = k\}$.
We can define a compact topology on the set $K^{E}$ of such functions. For this, endow $K$ with the co-finite topology where a set is closed iff it is finite or the whole of $K$.
Then endow $K^{E}$ with the product topology.

By \autoref{meetbond} a set $S$ is spanning for a matroid $M$ iff it meets every bond of that matroid. So we would like a function $f$ contained in each of the sets $C_{k,B} = \{f | S^f_k \cap B \neq \emptyset\}$, where $B$ is a bond for the matroid $M_k$. We will prove this by a compactness argument: we need to show that each $C_{k,B}$ is closed in the topology given above and that any finite intersection of them is nonempty.

To show that $C_{k,B}$ is closed, we rewrite it as $\bigcup_{e \in B}\{f | f(e) = k\}$. Each of the sets $\{f | f(e) = k\}$ is closed since their complements are basic open sets, and the union is finite since $M_k$ is co-finitary.

Now let $(k_i | 1 \leq i \leq n)$ and $(B_i| 1 \leq i \leq n)$ be finite families with each $B_i$ a bond in $M_{k_i}$. We need to show that $\bigcap_{1 \leq i \leq n} C_{k_i, B_i}$ is nonempty. Let $X = \bigcup_{1 \leq i \leq n} B_i$. Since the rank formula holds for each subset of $X$, we have by the finite version of the base packing Theorem a base packing $(S_k | k \in K)$ of $(M_k.X|k \in K)$. Now any $f$ such that $f(e) = k$ for $k \in S_k$ will be in $\bigcap_{1 \leq i \leq n} C_{k_i, B_i}$. This completes the proof.

\end{proof}

\autoref{pack_thm} does not extend verbatim to arbitrary infinite matroids.
Indeed, for every integer $k$ there exists a finitary matroid $M$ on a ground set $E$ with no three
disjoint bases yet satisfying $|Y | \geq kr_{M.Y}(Y)$ for every finite $Y \subseteq E$ \cite{aharoniThom, DiestelBook10}.

In the following we conjecture an extension of the finite base packing theorem to arbitrary infinite matroids.
This extension uses the following condition, which for finite sets $Y$ is implied by the rank formula of the base packing theorem but is still meaningful for infinite sets:

\begin{equation}\label{pack_cond}
\begin{minipage}[c]{0.8\textwidth}
If $(M_k.Y | k\in K)$ has a covering,
then it also has a packing.
\end{minipage}
\end{equation}

Indeed, if $(M_k.Y | k\in K)$ has a covering and there is an element of $Y$ contained in several of the corresponding independent sets, then this violates the rank formula.

Using our new condition, we obtain the following:

\begin{con}[Packing Conjecture]\label{pack_conj}
A family of matroids $(M_k|k\in K)$ on the same ground set $E$ has a packing
if and only if (\ref{pack_cond}) is true for every $Y\subseteq E$.
\end{con}

By a proof similar to that of \autoref{pc_is_cov_conj}, we obtain the following:

\begin{prop}\label{pc_is_pack_conj}
\autoref{pc} and \autoref{cov_conj} are equivalent. 
\end{prop}

As Packing/Covering is true for finite matroids, \autoref{pc_is_pack_conj} implies the non-trivial direction of \autoref{pack_thm}.
By \autoref{pack_cov_thm} we obtain the following applications.

\begin{cor}
A family of matroids $(M_k|k\in K)$ as in \autoref{pack_cov_thm} has a packing if and only if (\ref{pack_cond}) is true for every $Y\subseteq E$.
\end{cor}

In particular, we obtain the following:

\begin{cor}[Base packing theorem for the finite cycle matroid]
Any family of graphs $(G_k|k\in K)$ with a common edge set $E$ has a tree-packing if and only if 
(\ref{pack_cond_fc}) is true for every $Y\subseteq E$.
\begin{equation}\label{pack_cond_fc}
\begin{minipage}[c]{0.8\textwidth}
If $(M_k.Y | k\in K)$ has a tree-covering,
then it also has a tree-packing.
\end{minipage}
\end{equation}
\end{cor}

By \autoref{int_nf}, we also obtain the following.

\begin{cor}[Base packing theorem for the finite cycle matroid]
Any finite family of graphs $(G_k|k\in K)$ with edge set $E$ has a tree-packing if and only if (\ref{pack_cond_fc}) is true for every $Y\subseteq E$.
\end{cor}

A similar result was obtained by Aharoni and Ziv \cite{Aharoni:Ziv:92}.
However, their argument is different and they have the additional assumption that the ground set is countable.

Note that the covering conjecture for arbitrary finitary families is still open and equivalent to \autoref{o_fini}.

\section{Overview}
We have shown that a great many natural conjectures are equivalent, which we review here. The following are all equivalent.

\begin{description}
\item[The Intersection conjecture:]Any pair of matroids on the same ground set satisfies intersection
\item[The pairwise Packing/Covering conjecture:]Any pair of matroids on the same ground set satisfies Packing/Covering
\item[The Packing/Covering conjecture:]Any family of matroids on the same ground set satisfies Packing/Covering
\item[The Packing conjecture:]A family of matroids $(M_k|k\in K)$ on the same ground set $E$ has a packing
if and only if the following condition is true for every $Y\subseteq E$:
\begin{equation*}
\begin{minipage}[c]{0.8\textwidth}
If $(M_k.Y | k\in K)$ has a covering,
then it also has a packing.
\end{minipage}
\end{equation*}
\item[The Covering conjecture:]A family of matroids $(M_k|k\in K)$ on the same ground set $E$ has a covering
if and only if the following condition is true for every $Y\subseteq E$:
\begin{equation*}
\begin{minipage}[c]{0.8\textwidth}
If $(M_k\restric_Y | k\in K)$ has a packing,
then it also has a covering.
\end{minipage}
\end{equation*}
\end{description}

These equivalences allow the transfer of partial results (such as our proof of a special case of the Packing/Covering conjecture) to new contexts, and we hope that they will suggest new avenues for determining in what cases each of these conjectures holds.

\bibliographystyle{plain}
\bibliography{literatur}

\end{document}